\documentclass{amsart}

\usepackage{amsmath,amssymb,amsfonts,amstext,amsthm}

\newcommand{\nwc}{\newcommand}

\nwc{\aaa}{\mathcal{F}}
\nwc{\aap}{\mathcal{F}_{P}}

\nwc{\cb}{\overline{C}}
\nwc{\ccc}{\mathfrak{c}}
\nwc{\ch}{\widehat{C}}
\nwc{\cin}{\textbf{(v)}}
\nwc{\cl}{C'}
\nwc{\cp}{\mathcal{C}_{P}}
\nwc{\cpll}{\mathfrak{c}_{P'}}
\nwc{\ct}{\widetilde{C}}

\nwc{\dd}{\text{L}}
\nwc{\ddd}{\mathfrak{d}}
\nwc{\ddl}{\mathcal{L}'}
\nwc{\dlp}{\delta_{P}}
\nwc{\doi}{\textbf{(ii)}}

\nwc{\gtl}{\widetilde{g}}

\nwc{\hua}{h^{1}(C,\aaa )}

\nwc{\llb}{\mathcal{L}}

\nwc{\mm}{\mathfrak{m}}
\nwc{\mmp}{\mathfrak{m}_{P}}
\nwc{\mpd}{\mathfrak{m}_{P}^{2}}

\nwc{\nn}{\mathbb{N}}

\nwc{\ob}{\overline{\mathcal{O}}}
\nwc{\obr}{\mathcal{O}^*}
\nwc{\obp}{\overline{\mathcal{O}}_P}
\nwc{\och}{\mathcal{O}_{\hat{C}}}
\nwc{\oh}{\hat{\mathcal{O}}}
\nwc{\ohp}{\hat{\mathcal{O}}_{P}}
\nwc{\ol}{\mathcal{O}'}
\nwc{\oma}{\Omega (\mathfrak{a})}
\nwc{\omo}{\Omega (\mathcal{O})}
\nwc{\oo}{\mathcal{O}}
\nwc{\op}{\mathcal{O}_P}
\nwc{\opc}{\mathcal{O}_{P,C}}
\nwc{\oph}{\hat{\mathcal{O}}_{P}}
\nwc{\opl}{\mathcal{O}_{P}'}
\nwc{\oplc}{\mathcal{O}_{P,C}'}
\nwc{\opll}{\mathcal{O}_{P'}}
\nwc{\opt}{\tilde{\mathcal{O}}_{P}}
\nwc{\optt}{{\mathcal{O}}_{\tilde{P}}}
\nwc{\oq}{\mathcal{O}_{Q}}
\nwc{\oqt}{\tilde{\mathcal{O}}_{Q}}
\nwc{\ot}{\widetilde{\mathcal{O}}}
\nwc{\overop}{\bar{\oo}_{P}}

\nwc{\pb}{\overline{P}}
\nwc{\pbb}{P^*}
\nwc{\pbi}{\overline{P_{i}}}
\nwc{\pbr}{\overline{P_{r}}}
\nwc{\pgmd}{\mathbb{P}^{g+2}}
\nwc{\pgmu}{\mathbb{P}^{g+1}}
\nwc{\ph}{\hat{P}}
\nwc{\pp}{\mathbb{P}}
\nwc{\prv}{\noindent\textbf{Proof}:}
\nwc{\pt}{\widetilde{P}}
\nwc{\ptl}{\tilde{P}}
\nwc{\pum}{\mathbb{P}^{1}}

\nwc{\qh}{\hat{Q}}
\nwc{\qtl}{\tilde{Q}}
\nwc{\qua}{\textbf{(iv)}}

\nwc{\rh}{\hat{R}}

\nwc{\sei}{\textbf{(vi)}}
\nwc{\sep}{\beq\ast\ \ast\ \ast\enq}
\nwc{\ssp}{S_{P}}

\nwc{\tre}{\textbf{(iii)}}

\nwc{\um}{\textbf{(i)}}

\nwc{\vpb}{v_{\overline{P}}}
\nwc{\vtxp}{\widetilde{V}_{x,P}}
\nwc{\vxp}{V_{x,P}}

\nwc{\wh}{\hat{\omega}}
\nwc{\whp}{\hat{\omega}_{P}}
\nwc{\woch}{\omega\cdot\mathcal{O}_{\hat{C}}}
\nwc{\woh}{\omega\cdot\hat{\mathcal{O}}}
\nwc{\ww}{\omega}
\nwc{\wwb}{\omega^*}
\nwc{\wwct}{\omega _{\widetilde{C}}}
\nwc{\wwh}{\widehat{\omega}}
\nwc{\wwhp}{\widehat{\omega}_P}
\nwc{\wwp}{\omega _{P}}
\nwc{\wwt}{\widetilde{\omega}}
\nwc{\wwtp}{\widetilde{\omega}_P}

\nwc{\zz}{\mathbb{Z}}

\newtheorem{coro}{Corollary}[section]
\newtheorem{dfn}[coro]{Definition}

\newtheorem{lemma}[coro]{Lemma}

\newtheorem{rem}[coro]{Remark}
\newtheorem{thm}[coro]{Theorem}

\begin{document}

\title{Generalization of a Max Noether's Theorem}

\author{Renato Vidal Martins}
\address{Departamento de Matem\'atica, ICEx, UFMG
Av. Ant\^onio Carlos 6627,
30123-970 Belo Horizonte MG, Brazil}
\email{renato@mat.ufmg.br}



\begin{abstract}
Max Noether's Theorem asserts that if $\ww$ is the dualizing sheaf of a nonsingular nonhyperelliptic projective curve then the natural morphisms $\text{Sym}^nH^0(\omega)\to H^0(\omega^n)$ are surjective for all $n\geq 1$. This is true for Gorenstein nonhyperelliptic curves as well. We prove this remains true for nearly Gorenstein curves and for all integral nonhyperelliptic curves whose non-Gorenstein points are unibranch. The results are independent and have different proofs. The first one is extrinsic, the second intrinsic. 
\end{abstract}

\maketitle

\section{Introduction}

Let $C$ be an integral and complete curve of arithmetic genus $g$ over an algebraically closed field $k$. Let $\ww$ be its dualizing sheaf. We start by introducing the celebrated Max Noether's Theorem exactly as it is in \cite[pg 117]{ACGH}: if $C$ is nonsingular and nonhyperelliptic then the homomorphisms
$$
\text{Sym}^n\ H^0(C,\ww )\longrightarrow H^0(C,\ww ^n)
$$
are surjective for $n\geq 1$.

The same reference says the result is a consequence of projective normality of extremal curves, i.e., curves whose genus matches Castelnuovo's bound. Indeed, extremal curves are always projectively normal and this is a general fact proved in \cite[pp\,113-117]{ACGH} for nonsingular curves, but the reader should note the same proof holds for all integral curves as well. 

Now assume $C$ is Gorenstein. Then $\ww$ defines a morphism $\kappa :C\rightarrow\pp^{g-1}$. Let $C':=\kappa (C)$ be the canonical model of $C$. M. Rosenlicht proved in \cite{R} that $C'$ is extremal and that $\kappa$ is an isomorphism if $C$ is nonhyperelliptic. Therefore Max Noether's Theorem holds actually for all Gorenstein nonhyperelliptic curves. 

We start this article with Theorem \ref{thmext} which proves Max Noether's result in the same way as above but for a bigger bunch of curves. The first difference is we derive projective normality from linear normality. Most of the work was already done in \cite[Lem.\,5.4.(1)]{KM} which we rephrase and reprove in Lemma \ref{lemprj}. The second difference is, in order to extend the result, we deal with $\widehat{C}$, the blowup along $\ww$, instead of $C$. The new curves which appear were called nearly Gorenstein in \cite{KM}.  They have just one non-Gorenstein point and for which the local ring is almost Gorenstein, a desirable (though restrictive) property introduced by R. Fr\"oberg and V. Barucci in \cite{BF}. 

We think that Theorem \ref{thmext} is the best one gets with an extrinsic argument. But we recall that, first of all, Max Noether's Theorem has to do with the dualizing sheaf itself no matter where it embeds the curve or its blowup.  So in order to deal with this problem intrinsically we were motivated by a strong result due to M. Rosenlicht \cite[Thm.\,17]{R}. It asserts that, if $C$ is nonhyperelliptic, then the birational map between $C$ and the canonical model $C'$ [cf.\,Definition\,\ref{dfncan}] is regular on $C'$. His proof, specially in the non-Gorenstein case, is focused on the local rings of this sort of points and the stalks of the dualizing sheaf at them; and the technique was the computation of values of differentials. Indeed, in the last two pages of the article, where is the core of Rosenlicht's proof, there's only one paragraph where he does not compute values. So with this same tool we prove the following statement, which is our main result: if $C$ is a nonhyperelliptic curve whose non-Gorenstein points are unibranch then Max Noether's assertion holds (Theorem \ref{thmmon}).

The hypothesis assumed that the non-Gorenstein points are unibranch is due to property (\ref{equvwp}) below which describes the possible values of differentials that are regular at a given point of $C$. In the general case, it is not easy to deal with this property \cite[Thm.\,2.11]{S} -- as Rosenlicht did -- within our context. In fact, we point out that Max Noether's assertion is stronger than Rosenlicht's one [cf. Remark \ref{remext}]. In other words, our task is a little harder. This led us to think that the multibranch problem deserves another work. In this one, the reader should note that, although in Section \ref{secint} we always assume the non-Gorenstein points are unibranch, it is just Step 2 of the proof of Lemma \ref{lemtch} what really must be extended in order to get a sharper result. But the question if Max Noether's statement holds for all integral curves remains open for us.

\ 

\paragraph{\bf Acknowledgments.}
This article continues \cite{KM}, a joint work  with Steven L. Kleiman, to whom we thank very much for an invitation to MIT and regular email discussion after that when many suggestions were built in to this paper. The author is partially supported by CNPq grant number PDE 200999/2005-2.

\

\section{Generalization with an extrinsic argument}

Let $C$ be an integral and complete curve of arithmetic genus $g$ over an algebraically closed field $k$ of arbitrary characteristic. Let $\oo :=\oo _C$ be the structure sheaf on $C$. According to \cite[pg.\,140]{ACGH} we define.

\begin{dfn}
\emph{A curve $C\subset\pp^r$ is $n$-\emph{normal} if the hypersurfaces of degree $n$ cut out the complete linear series $|\oo_C(n)|$. A $1$-normal curve is said \emph{linearly normal} and a curve is said \emph{projectively normal} if it is $n$-normal for all $n\geq 1$.}
\end{dfn}

\begin{lemma}
\label{lemprj}
Let $C$ be a nondegenerate curve of degree $d$ in $\pp^r$. If $d<2r$ or else $d=2r$ and $h^1(\oo_C(1))>0$ then $C$ is linearly normal iff it is projectively normal. 
\end{lemma}

\begin{proof}
Sufficiency is immediate. To prove necessity, consider the left exact sequence
\begin{equation}
\label{eqLP3}
0\longrightarrow H^0(\oo_{C}(n-1))\stackrel{u}\longrightarrow H^0(\oo_C(n))\stackrel{v}\longrightarrow H^0(\oo_H(n)), 
\end{equation}
where $H$ is a hyperplane divisor on $C$.

Let $V_n$ denote the image of $H^0(\oo_{\pp^{r}}(n))$ in $H^0(\oo_C(n))$, and set $W_n:=v(V_n)$.  Then $\dim(W_n)\geq \min\{d,\ n(r-1)+1\}$ by \cite[Lem.,\,pg.\,115]{ACGH}.  Hence
\begin{equation}
\label{eqLP1}
h^0(\oo_C(n))-h^0(\oo_C(n-1))\geq\min\{d,\ n(r-1)+1\}.
\end{equation}

Also, if equality holds in (\ref{eqLP1}), then $v(H^0(\oo_C(n)))=W_n$ since both sides have the same dimension.  So
$H^0(\oo_C(n))$ is spanned by $V_n$ and $\text{Im}(u)$.  But
$u(V_{n-1})\subset V_n$.  And, if $C$ is linearly normal, then
$V_1=H^0(\oo_C(1))$.  Hence, if in addition, equality holds in
(\ref{eqLP1}) for $n\ge2$, then induction on $n$ yields $V_n=H^0(\oo_C(n))$ for $n\geq 1$; in other words, then $C$ is projectively
normal.  Thus to complete the proof, we have to prove that
equality holds in (\ref{eqLP1}) for $n\ge2$.

Set $h(n):=h^1(\oo_C(n))$.  Then $h^0(\oo_C(n))=nd+1-g+h(n)$
by the Riemann--Roch Theorem.  Hence the bound~(\ref{eqLP1}) is
equivalent to this bound:
\begin{equation}
\label{eqLP2}
d-(h(n-1)-h(n))\ge\min\{d,\ n(r-1)+1\}.
\end{equation}
Here $h(n-1)-h(n)\ge0$ because the sequence~(\ref{eqLP3}) continues,
ending with $$H^1(\oo_C(n-1))\to H^1(\oo_C(n))\to 0.$$
Write $d=2r-a$ for $a\geq 0$. We have
\begin{equation*}
\label{eqeqLP4}
(n(r-1)+1)-d=(n-2)(r-1)-1+a.
\end{equation*}
If $r=1$ then $C=\pp^1$ which is projectively normal. So assume $r\geq 2$. For $n\ge2$, the right side is nonnegative unless $a=0$ and $n=2$. 

Hence, for $n\geq 3$, the right side of (\ref{eqLP2}) is equal to $d$.
But $h(n-1)-h(n)\geq 0$.  Therefore, equality holds in (\ref{eqLP2}), and
$h(n-1)=h(n)$.  But, by Serre's Theorem,
$h(n)=0$ for $n\gg0$.  So $h(n)=0$ for $n\ge2$.

Suppose $a>0$.  Then similarly, equality holds in (\ref{eqLP2}) for
$n=2$ too, and $h(1)=0$.  So equality holds in (\ref{eqLP1}) for
$n\ge2$, as desired.  

Finally, instead suppose $a=0$. Then $d=2r$ and $h(1)>0$ by hypothesis. Now, take $n=2$ in (\ref{eqLP2}), getting $2r-h(1)$ on
the left as $h(2)=0$, and $2r-1$ on the right.  Hence, $h(1)=1$ and
equality holds in (\ref{eqLP2}) for $n=2$.  So equality holds in
(\ref{eqLP1}) for $n\ge2$, as desired. The proof is now complete.
\end{proof}

A curve $C$ is said \emph{hyperelliptic} if there is a morphism $C\rightarrow\pp^1$ of degree $2$. Let $\ww_C$, or simply $\ww$, denote the dualizing sheaf. A curve $C$ is said \emph{Gorenstein} if $\ww$ is invertible. A point $P\in C$ is said \emph{Gorenstein} if $\ww_P$ is a free $\oo_P$-module.

Given any integral scheme $A$, any map $\alpha :A\to C$ and a sheaf $\mathcal{G}$ on $C$, set
$$\mathcal{O}_A\mathcal{G}:= \alpha^*\mathcal G/\text{Torsion}(\alpha^*\mathcal G).$$

Let $\nu :\overline{C}\rightarrow C$ be the normalization map. Set $\overline{\oo}:=\nu _{*}(\oo _{\overline{C}})$. We denote $\mathcal{C}$ the conductor of $\overline{\oo}$ into $\oo$. We also set $\overline{\ww}:=\nu_*(\ww_{\overline{C}})$ and $\overline{\oo}\ww:=\nu_*(\oo_{\overline{C}}\ww$).

Given any coherent sheaf $\mathcal{F}$ on $C$ set $\mathcal{F}^n:=\text{Sym}^n\mathcal{F}/\text{Torsion}(\text{Sym}^n\mathcal{F})$. If $\mathcal{F}$ is invertible then clearly $\mathcal{F}^{n}=\mathcal{F}^{\otimes n}$. Let $\phi :\text{Sym}^n\,\mathcal{F}\to H^0(\mathcal{F}^n)$ be the natural morphism. We set $H^0(\mathcal{F})^n:=\phi(\text{Sym}^n\,\mathcal{F})$.

Call $\widehat{C}:=\text{Proj}(\oplus\,\ww ^n)$ the blowup of $C$ along $\ww$. Let $\beta :\widehat{C}\rightarrow C$ be the natural morphism. Set $\widehat{\oo}=\beta_*(\oo_{\widehat{C}})$ and $\widehat{\oo}\ww:=\beta_*(\oo _{\widehat{C}}\ww)$.

\begin{dfn}
\label{defvxp}
\emph{For $E\subset k(C)$ set $E^n:=\{\sum_{j\in J} f_{1j}\ldots f_{nj}\,|\,f_{ij}\in E\ \text{and}\ |J|<\infty\}$.
Let $z\in\Omega_C$ be a differential. Set 
$$
\begin{array}{lll}
W_z:=H^0(\ww)/z\ \ \ & V_{z,P}:=\wwp/z\ \ \ & W_{z,n}:=\displaystyle{\bigcap_{P\in C}}\,V_{z,P}^n.
\end{array}
$$
If $\varphi :\ct\to C$ is a morphism and $\wwt:=\varphi_{*}(\ww_{\ct})$ set}
$$
\begin{array}{lll}
\widetilde{W}_z:=H^0(\wwt)/z\ \ \ & \widetilde{V}_{z,P}:=\wwtp/z\ \ \ & \widetilde{W}_{z,n}=\displaystyle{\bigcap_{P\in C}}\,\widetilde{V}_{z,P}^n.
\end{array}
$$
\end{dfn}

\begin{dfn}
\label{deffix}
\emph{We fix, throughout this paper, a differential $x\in H^0(\ww)$ such that $(\overline{\oo}\ww)_P=\overline{\oo}_Px$ for every singular point $P\in C$. Such a differential exists because $H^0(\ww)$ generates $\overline{\oo}\ww$ as proved in \cite[p.\,188 top]{R}, the singular points of $C$ are of finite number and $k$ is infinite since it is algebraically closed.}
\end{dfn} 

\begin{dfn}
\label{dfnlyG}
\emph{Call $C$ {\it nearly Gorenstein\/} if $C$ has only one non-Gorenstein point $P$ and if the local ring $\oo_P$
  is almost Gorenstein in the sense of Barucci and Fr\"oberg
  \cite[p.\,418]{BF}, namely, if
$$\dim(\overline{\oo}_P/\oo_P) = \dim (\oo_P/\mathcal{C}_P)+\dim(\text{Ext}^1(k,
\oo_P))-1$$
where $k$ is the ground field.}
\end{dfn} 

\begin{thm}
\label{thmext}
If $C$ is either nonhyperelliptic Gorenstein or else nearly Gorenstein then the homomorphisms
$$
\text{\emph{Sym}}^n\,H^0(C,\ww )\longrightarrow H^0(C,\ww ^n)
$$
are surjective for $n\geq 1$.
\end{thm}

\begin{proof}
From \cite[Prp.\,4.5]{KM} we have that $\oo_{\widehat{C}}\ww$ is an ample invertible sheaf on $\widehat{C}$ which is generated by $H^0(\ww)$. So, first, $\widehat{C}=\text{Proj}(\oplus\,H^0((\oo_{\widehat{C}}\ww)^{\otimes n}))$ and, second, the complete linear system $\text{L}:=|\oo_{\widehat{C}}\ww|$ is base point point free since $H^0(\ww)\subset H^0(\oo_{\widehat{C}}\ww)$. Set $r:=h^0(\oo_{\widehat{C}}\ww)-1$ and let $\varphi :\widehat{C}\rightarrow\pp^{r}$ be the morphism defined by $\text{L}$. Call $C^*:=\varphi(\widehat{C})$. Now $C^*=\text{Proj}(\oplus\,H^0(\oo_{C^*}(n)))$ by the very definition of $C^*$. Hence, for every $n\geq 1$, we have a sequence of linear morphisms
\begin{equation}
\label{equscw}
\text{Sym}^n\ H^0(\oo _{\widehat{C}}\ww)\stackrel{\alpha_{n}}\longrightarrow H^0(\oo _{C^*}(n))\stackrel{\beta_{n}}\longrightarrow H^0((\oo_{\widehat{C}}\ww)^{\otimes n})
\end{equation}
where $\alpha _n$ is the natural homomorphism and $\beta_n$ is the injective homomorphism determined by $\varphi$. 

If $C$ is Gorenstein then the $\alpha_n$'s are surjective. In fact, set $d:=\text{deg}(C^*)$. If $C$ is Gorenstein, then $\widehat{C}=C$; $\oo_{\widehat{C}}\ww=\ww$; $r=g-1$; the morphism $\varphi$ agrees with the canonical morphism $\kappa: C\rightarrow \pp^{g-1}$ and $C^*=C'$, the canonical model. Then either $d<2r$ and hence $C^*$ is projectively normal by Lemma \ref{lemprj} and hence the $\alpha_n$'s are surjective; or else,
$d=2r=2g-2$. In this case, $\kappa$ is an isomorphism owing to \cite[Thm.\,4.3]{KM}. Hence $\oo_{C^*}(1)=\ww$ and thus $h^1(\oo_{C^*}(1))=1$. Lemma \ref{lemprj} now implies $C^*$ is projectively normal and the $\alpha_n$'s are surjective.

If $C$ is Gorenstein and nonhyperelliptic then the $\beta_n$'s are surjective. In fact, if $C$ is Gorenstein and nonhyperelliptic then $\varphi$ is an isomorphism owing to \cite[Thm.\,4.3]{KM}. Hence the $\beta_n$'s are surjective.

So the result is proved for Gorenstein hyperelliptic curves. Let us prove it for nearly Gorenstein curves.

If $C$ is non-Gorenstein then the $\alpha_n$'s and the $\beta_n$'s are surjective. In fact, if $C$ is non-Gorenstein then $\oo_{\widehat{C}}\ww$ is very ample and $h^1(\oo_{\widehat{C}}\ww)=0$ owing to \cite[Prp.\,5.2]{KM}. Besides, from the proof of \cite[Prp.\,5.2]{KM}, we also have $\widehat{g}\leq g-2$ where $\widehat{g}$ is the genus of $\widehat{C}$. Hence, first, $\varphi$ is an isomorphism and the $\beta_n$'s are surjective. Second,
\begin{align*}
d & =\text{deg}_{\widehat{C}}(\oo_{\widehat{C}}\ww)=h^0(\oo_{\widehat{C}}\ww)-1+\widehat{g}     \\             
  & =r+\widehat{g}\leq r+g-2=r+h^0(\ww)-2                                         \\
  & \leq r+h^0(\oo_{\widehat{C}}\ww)-2                                            \\
  & =r+(r+1)-2<2r.
\end{align*} 
Therefore Lemma \ref{lemprj} implies $C^*$ is projectively normal, that is, the $\alpha_n$'s are surjective.

Now assume $C$ is nearly Gorenstein. From \cite[Lem.\,5.8]{KM} we have $H^0(\oo_{\widehat{C}}\ww)=H^0(\ww)$. So it suffices to prove that  
$H^0((\oo_{\widehat{C}}\ww)^{\otimes n})=H^0(\ww^n)$ for $n\geq 2$. But we claim that, in this case, actually $(\widehat{\oo}\ww)^{\otimes n}=\ww^n$ for $n\geq 2$. In fact, if $C$ is nearly Gorenstein then $C$ has just one non-Gorenstein point $P$, and for which the local ring $\oo_P$ is almost Gorenstein. Moreover, from \cite[Prp.\,28]{BF} we have that if $\op$ is almost Gorenstein then $\text{dim}((\widehat{\oo}\ww)_{P}/\ww_P)=1$. Now, for the differential $x$ of Definition \ref{deffix}, we have $(\widehat{\oo}\ww)_P=\widehat{\oo}_Px$. Besides, $H^0(\ww)$ generates $\ww$ owing to \cite[p.\,536\,mid]{EKSH}. This implies $\op\subset\vxp\subset\widehat{\oo}_P$. Then $\vxp^2\subset\widehat{\oo}_P$ owing to \cite[Lem\,6.1.(b)]{KM} and $\vxp\subsetneq\vxp^2$ owing to the proof of \cite[Prp.\,28]{BF}. Hence $\vxp^2=\widehat{\oo}_P$ because $\text{dim}(\widehat{\oo}_{P}/\vxp)=\text{dim}((\widehat{\oo}\ww)_{P}/\ww_P)=1$. Therefore $\vxp^n=\widehat{\oo}_P$ for $n\geq 2$. This yields $\wwp^n=(\widehat{\oo}\ww)^{\otimes n}_P$ for $n\geq 2$. Since $P$ is the only non-Gorenstein point of $C$, then $\ww$ and $\widehat{\oo}\ww$ agree outside $P$. Thus $(\widehat{\oo}\ww)^{\otimes n}=\ww^n$ for $n\geq 2$ and the claim is proved. Now $H^0((\widehat{\oo}\ww)^{\otimes n})=H^0((\oo_{\widehat{C}}\ww)^{\otimes n})$ and the theorem is proved.
\end{proof}

\begin{dfn}
\label{dfncan}
\emph{In \cite[p.\,188\,top]{R} Rosenlicht showed that the linear system $\overline{\text{L}}:=(\oo_{\overline{C}}\ww,H^0(\ww))$ is base point free. He considered then the morphism $\overline{\kappa}:\overline{C}\rightarrow\pp^{g-1}$ defined by $\overline{\text{L}}$ and called $C':=\overline{\kappa}(C)$ the \emph{canonical model} of $C$. In \cite[Dfn.\,4.9]{KM} one finds another characterization of $C'$. It is the image of the morphism $\widehat{\kappa}:\widehat{C}\rightarrow\pp^{g-1}$ defined by the linear system $\widehat{\text{L}}:=(\oo_{\widehat{C}}\ww,H^0(\ww))$. Note this definition of $C'$ agrees with the former, which appears in the Introduction. In fact, if $C$ is Gorenstein then $\widehat{C}=C$, $\oo_{\widehat{C}}\ww=\ww$ and $\widehat{\kappa}$ is nothing but $\kappa$, the canonical morphism.}
\end{dfn}

\begin{rem}
\label{remext}
\emph{As pointed out in the Introduction, Max Noether's assertion is stronger than Rosenlicht's one. In fact, since $\ww$ is generated by global sections, Rosenlicht's assertion is equivalent to this one: $\widehat{\kappa}:\widehat{C}\rightarrow C'$ is an isomorphism. So assume $C$ is nonhyperelliptic and assume Max Noether's assertion holds, that is, $\text{Sym}^n\ H^0(\ww )\to H^0(\ww ^n)$ is surjective for $n\geq 1$. We will show $\widehat{\kappa}$ is an isomorphism. We first claim there exists $n$ such that $\ww^n =(\widehat{\oo}\ww)^{\otimes n}$. Indeed, for every singular point $P\in C$ holds $\oo_P\subset\vxp\subset\widehat{\oo}_{P}$. Consider the sequence
$$
\vxp\subset\vxp^2\subset\ldots\subset\vxp^i\subset\ldots\subset\widehat{\oo}_P.
$$
There exists $n$, depending on $P$, such that $\vxp^i=\vxp^n$ for every $i\geq n$ because $\text{dim}(\widehat{\oo}_{P}/\vxp)<\infty$. Now $\vxp^n$ is a ring and $\widehat{\oo}_P$ is the smallest ring in $k(C)$ which contains $\vxp$ due to \cite[Lem.\,6.1.(b)]{KM}. Thus $\vxp^n$ and $\widehat{\oo}_P$ agree and so do $\wwp^n$ and $(\widehat{\oo}\ww)^{\otimes n}$. Take $n$ which works for all singular points of $C$. Since $\ww$ and $\widehat{\oo}\ww$ agree outside the singular points of $C$ we have $(\widehat{\oo}\ww)^{\otimes n}=\ww^n$. For this $n$, consider the sequence
$$
\text{Sym}^n\ H^0(\ww)\stackrel{\alpha}\longrightarrow H^0(\oo _{C'}(n))\stackrel{\beta}\longrightarrow H^0((\oo_{\widehat{C}}\ww)^{\otimes n})=H^0((\widehat{\oo}\ww)^{\otimes n})=H^0(\ww^n)
$$
where $\alpha$ and $\beta$ are naturally associated to the morphism $\widehat{\kappa}:\widehat{C}\rightarrow C'$. Max Noether's assertion implies $\beta\circ\alpha$ is surjective, so $\beta$ is surjective. But $\beta$ is injective, so $\beta$ is bijective. Hence $h^0(\oo_{C'}(n))=h^0((\widehat{\oo}\ww)^{\otimes n})$. If $C$ is nonhyperelliptic then $\deg _{C'}(\oo_{C'}(n))=\deg _{\widehat{C}}(\widehat{\oo}\ww)^{\otimes n}$. Besides, $h^1(\oo_{C'}(n))=h^1((\widehat{\oo}\ww)^{\otimes n})=0$ taking $n>>0$. Hence the Riemann-Roch theorem implies that $\widehat{C}$ and $C'$ are of the same arithmetic genus.  Therefore, $\widehat{\kappa}:\widehat{C}\to C'$ is an isomorphism; in other words, Rosenlicht's assertion holds.}
\end{rem}

\section{Generalization with an intrinsic argument (Main Theorem)}
\label{secint}

In this section we show our Main Theorem announced in the Introduction, that is, Max Noether's assertion holds for nonhyperelliptic curves for which the non-Gorenstein singularities satisfie the following property. 

\begin{dfn}
\emph{A point $P\in C$ is called \emph{unibranch} if $\pi^{-1}(P)$ consists of only one point.}
\end{dfn}

If $P\in C$ is unibranch, fix a local parameter $t_P$ of the local ring $\obp$. We define $\alpha_P ,\beta_P\in\nn$ such that 
\begin{equation}
\label{equalb}
\cp =\obp t_P^{\alpha_P}\ \ \ \text{and}\ \ \ \obp\mmp =\obp t_P^{\beta_P}
\end{equation}
where $\mmp$ is the maximal ideal of $\op$. Note $\beta_P$ agrees with the multiplicity of $P$. We use the same notation $v_P$ to the valuation of $\obp$ applied to either rational functions or differentials.

\begin{lemma}
\label{lemtch}
Let $P\in C$ be a unibranch non-Gorenstein point. Then 
$$
W_x^n\longrightarrow\ \vxp^n\,/\,t_P^{-\epsilon}\cp^n
$$
is surjective for every $n\geq 1$, where
\begin{enumerate}
\item[(i)] $\epsilon=2n-1$
\item[(ii)] $\epsilon=1$ if there is $y_0\in H^0(\ww)$ with $v_P(y_0)=0$.
\item[(iii)] $\epsilon=0$ if there are $y_0,y_1\in H^0(\ww)$ with $v_P(y_0)=0$ and $v_P(y_1)=1$ or $2$.
\end{enumerate}
\end{lemma}

\begin{proof}
Call $\alpha :=\alpha_P$, $\beta :=\beta_P$ and $t:=t_P$ for short. Consider the sequence 
$$
\vxp^n\supset\cp\supset t^{\alpha-\beta}\cp =t^{-\beta}\cp ^2\supset\cp^2\supset\cp ^n
$$
of which each inclusion will correspond to a step of our proof.

\noindent \emph{Step 1}. $W_x^n\twoheadrightarrow\vxp^n/\cp$.

From the proof of \cite[Lem.\,6.1]{KM} we have $\vxp\subset W_x+\cp$. On the other hand, both $\cp$ and $W_x$ are contained in $\vxp$. We are led to
\begin{equation}
\label{equwhc}
\vxp=W_x+\cp. 
\end{equation}

Now $W_x\subset\vxp\subset\overline{\oo}_P$. Thus if $a\in W_x$ and $b\in\cp$ then $ab\in\cp$. This yields
\begin{equation}
\label{equwpl}
\vxp ^n= W_x^n+\cp 
\end{equation}
which proves the claim of first step.
 
\noindent \emph{Step 2}. $\cp /t^{\alpha-\beta}\cp$ is generated by images of elements in $W_x^2$.

Call $v:=v_P$ for short, set $\text{S}:=v(\op )$. From \cite[Thm.\,2.11]{S} we have
\begin{equation}
\label{equvwp}
v(\wwp )=\{ d\in\zz\ |-d-1\not\in \text{S} \}.  
\end{equation}
Now $\overline{\ww}_P=\mathcal{C}_Px$ owing to \cite[Lem.\,2.8]{KM} which implies $v(x)=-\alpha$. This yields 
\begin{equation}
\label{equvwl}
v(\vxp )=\text{K}:=\{ d\in\zz\ |\ \alpha-d-1\not\in \text{S}\}.  
\end{equation}
Since $P$ is non-Gorenstein, $v(\vxp )\supsetneq v(\op )$, that is, $\text{K}\supsetneq\text{S}$. Take $d_1\in \text{K}\setminus \text{S}$. Then $d_2:=\alpha -d_1-1\in \text{K}\setminus \text{S}$ by the very definition of $\text{K}$.

Let $r$ be the greatest integer such that $(r+1)\beta\leq\alpha$. Assume $r\geq 1$. In particular, $2\beta\leq\alpha$. We claim that for every $1\leq i\leq r$ one can find natural numbers $q_{i1},q_{i2}$ such that $i=q_{i1}+q_{i2}$ and $q_{ij}\beta +dj<\alpha$ for $j=1,2$. Indeed, it suffices to prove this for $r$ because if $i\leq r$ one can take $q_{i1}:=\text{min}\{ i,q_{r1}\}$ and $q_{i2}:=i-q_{i1}\leq q_{r2}$. Suppose without loss of generality $d_1\leq d_2$. Take $q_{r2}$ as the greatest integer such that $q_{r2}\beta\leq d_1$ and $q_{r1}:=r-q_{r2}$. Since $d_1\leq d_2$ and $d_1+d_2=\alpha -1$ we have $2d_1<\alpha$. Therefore $(q_{r2} +1)\beta\leq d_1+\beta\leq 2\,\text{max}\{ d_1,\beta\} \leq\alpha$ and hence $r\geq q_{r2}$ due to the definition of $r$. This implies $q_{r1}\geq 0$. Now $q_{r2}\beta +d_2\leq d_1+d_2<\alpha$ and $q_{r1}\beta +d_1=(r-q_{r2})\beta +d_1=r\beta -q_{r2}\beta +d_1\leq(\alpha -\beta )+ (d_1-q_{r2}\beta )< \alpha$ because $d_1-q_{r2}\beta<\beta$ due to the definition of $q_{r2}$. This proves the claim.

Now take $a_1,a_2\in\vxp$ such that $v(a_1)=d_1$ and $v(a_2)=d_2$, and $m\in\op$ such that $v(m)=\alpha$. Since $\vxp$ is an $\op$-module, $m^{q_{ij}}a_j\in\vxp$ for $1\leq i\leq r$ and $j=1,2$. From (\ref{equwhc}), write $m^{q_{ij}}a_j=a_{ij}'+a_{ij}''$ with $a_{ij}'\in W_x$ and $a_{ij}''\in\cp$. Set $f_i:=a_{i1}'a_{i2}'$. We have $f_i\in W_x^2$ for $1\leq i\leq r$. Besides, 
\begin{align*}
v(f_i ) &= v(a_{i1}')+v(a_{i2}')=v(m^{q_{i1}}a_1)+v(m^{q_{i2}}a_2)\\
        &= v(m^{q_{i1}+q_{i2}}a_1a_2)= v(m^ia_1a_2) \\
        &= i\beta +d_1+d_2=i\beta +\alpha -1
\end{align*}          
where the second equality holds because, for $j=1,2$, we have $v(m^{q_{ij}}a_j)=q_{ij}\beta +d_j<\alpha\leq v(a_{ij}'')$ what implies $v(m^{q_{ij}}a_j)=v(a_{ij}')$. 

Combining  (\ref{equwhc}) and (\ref{equvwl}) one can find a set of elements $\{b_1,\ldots ,b_{\beta -1}\}\subset W_x$ such that 
$$
v(b _i)=i+\alpha-\beta -1.
$$
Now for $1\leq i\leq r$ we clearly have $m^i\in\op\subset\vxp$ and we also have $v(m^i)=i\beta <\alpha$. Use (\ref{equwhc}) to find $m_i\in W_x$ such that $v(m_i)=i\beta$. If $\beta\not |\alpha$ one also finds $m_{r+1}\in W_x$ such that $v(m_{r+1})=(r+1)\beta$. Write $\alpha =(r+1)\beta +p$ with $0\leq p<\beta$. Set
$$
A_1:=\left(\bigcup_{i=1}^{r}\{ m_i b_1,\ldots ,m_i b_{\beta -1},f_i\}\right)\ \bigcup\ \{ m_{r+1}b_1,\ldots ,m_{r+1}b_p\}
$$
where the set at the left hand side of the union is empty if $r=0$ and so is the set at the right hand side if $p=0$. We claim $A_1$, which is contained in $W_x^2$ by construction, provides a basis for $\cp /t^{\alpha-\beta}\cp$. In fact, first, $A_1$ is of the right size because $\dim (\cp /t^{\alpha-\beta}\cp )=\alpha -\beta$ and $|A_1|=r\beta +p=\alpha -\beta$. Second, $A_1\subset\cp$ and its elements are linearly independent mod $t^{\alpha-\beta}\cp$ as one can see computing values:  
$$
\begin{array}{llll}
v(m_1b_1)=\alpha & \ldots & v(m_1b_{\beta -1})=\alpha+\beta-2  & v(f_1)=\alpha +\beta -1  \\
\ \ \ \ \ \ \ \vdots         &  &\ \ \ \ \ \ \ \vdots  & \ \ \ \ \ \ \ \vdots\\
v(m_rb_1)=\alpha +(r-1)\beta & \ldots & v(m_rb_{\beta -1})=\alpha+r\beta -2 & v(f_r)=\alpha+r\beta -1  \\ 
v(m_{r+1}b_1)=\alpha+r\beta & \ldots & v(m_{r+1}b_{p})=2\alpha -\beta -1. &
\end{array}
$$
This proves the claim and we are done with second step.

\noindent\emph{Step 3.1.} $t^{-\beta}\cp^2/t^{-3}\cp^2$ is generated by images of elements in $W_x^2$.

Assume $\beta>3$ for otherwise the claim triviously holds. Consider the set
$$
A_2 :=\{ b_{\beta-1}b_3,\ldots,b_{\beta-1}b_{\beta-2},b_{\beta-1}^2\}. 
$$
It is contained in $W_x^2$ by construction and provides a basis for $t^{-\beta}\cp ^2/t^{-3}\cp^2$. In fact, $\dim (t^{-\beta}\cp^2/t^{-3}\cp ^2)=\beta-3=|A_2|$. Besides, $A_2\subset t^{-\beta}\cp^2$ and its elements are linearly independent mod $t^{-3}\cp^2$ as one can see computing values:  
$$
\begin{array}{llll}
v(b_{\beta-1}b_3)=2\alpha -\beta  & \ldots &  v(b_{\beta-1}^2)=2\alpha -4.    
\end{array}
$$
\noindent \emph{Step 3.2}. If (ii) or (iii) holds, then $t^{-\beta}\cp^2/t^{-\epsilon}\cp^2$ is generated by images of elements in $W_x^2$.

If there is $y_0\in H^0(\ww)$ with $v(y_0)=0$, then $h_0:=y_0/x\in W_x$ and $v(h_0)=\alpha$. Consider the set
$$
A_2':=\{ h_0b_1,\ldots,h_0b_{\beta-1}\}.
$$
It is contained in $W_x^2$ by construction and provides a basis for $t^{-\beta}\cp ^2/t^{-1}\cp^2$. In fact, $\dim (t^{-\beta}\cp^2/t^{-1}\cp ^2)=\beta-1=|A_2'|$. Besides, $A_2\subset t^{-\beta}\cp^2$ and its elements are linearly independent mod $t^{-1}\cp^2$ as one can see computing values:  
$$
\begin{array}{lll}
v(h_0b_1)=2\alpha -\beta  & \ldots &  v(h_0b_{\beta-1})=2\alpha -2.    
\end{array}
$$
If, besides, there is $y_1\in H^0(\ww)$ such that $v(y_1)=1$ (resp. $v(y_1)=2$) then $h_1:=y_1/x\in W_x$ and $v(h_1)=\alpha+1$ (resp. $v(h_1)=\alpha+2$). Thus one forms $A_2''$ adjointing $h_1b_{\beta-1}$ (resp. $h_1b_{\beta-2}$) to $A_2'$ in order to get a basis for $t^{-\beta}\cp ^2/\cp^2$. 

\noindent\emph{Step 4.1.} $t^{-3}\cp^2/t^{-2n+1}\cp^n$ is generated by images of elements in $W_x^n$.

The set $B:=\{f_1\}\ \bigcup A_1\bigcup A_2$ is contained in $W_x^2$ by construction. Then
$$
A_3 :=\bigcup _{i=1}^{n-2}\ b_{\beta-1}^i\,B
$$ 
and $A_3\subset t^{-3}\cp^2$ and its elements are linearly independent mod $t^{-2n+1}\cp^n$ as one sees computing values:
$$
\begin{array}{lll}
v(b_{\beta-1}f_1)=2\alpha -3  & \ldots &  v(b_{\beta-1}^{n})=n\alpha-2n.
\end{array}
$$

\noindent\emph{Step 4.2.} If (ii) or (iii) holds, then $t^{-\epsilon}\cp^2/t^{-\epsilon}\cp^n$ is generated by images of elements in $W_x^n$.

The set $B':=\{f_1\}\ \bigcup A_1\bigcup A_2'$ is contained in $W_x^2$ by construction. Then
$$
A_3' :=\bigcup _{i=1}^{n-2}\ \ h_0^i\,B'
$$ 
is contained in $W_x^n$ and it is easily checked that $A_3\subset t^{-1}\cp^2$ and its elements are linearly independent mod $t^{-1}\cp^n$. If there is $y_1\in H^0(\ww)$ such that $v(y_1)=1$ or $2$, form $B''$ discarding $f_1$ in $B'$ and replacing $A_2'$ by $A_2''$. Then form $A_3''$ replacing $B'$ by $B''$ in $A_3'$. We are finally done.
\end{proof}

\begin{lemma}
\label{lemvpw}
Let $P\in C$ be a nonsingular point. Then
\begin{enumerate}
\item[(i)] If $g\geq 1$, there is $y_0\in H^0(\ww)$ with $v_P(y_0)=0$.
\item[(ii)] If $C$ is nonhyperelliptic, there is $y_1\in H^0(\ww)$ with $v_P(y_1)=1$.
\item[(iii)] If $C$ is hyperelliptic and $g\geq 2$, there is $y_1\in H^0(\ww)$ with $v_P(y_1)=1$ or $2$.
\end{enumerate}
\end{lemma}

\begin{proof}
If $g\geq 1$ then $\ww$ is generated by global sections, in particular $\ww_P=\op y_0$ for $y_0\in H^0(\ww)$. Since $P$ is nonsingular, from (\ref{equvwp}) there is $y\in\wwp$ such that $v_P(y)=0$. We have $v_P(y_0)\geq 0$ due to (\ref{equvwp}) as well. Since $y=fy_0$ for $f\in\op$, this forces $v_P(y_0)=0$ and (i) is proved. 

Assume $C$ is nonhyperliptic. If $h^0(\ww(-P))=h^0(\ww(-2P))$ then $h^1(\ww(-2P))>h^1(\ww(-P))$ because $\deg (\ww(-P))>\deg(\ww(-2P))$. So $h^0(\oo(2P))>h^0(\oo(P))\geq 1$. Since $\deg(\oo(2P))=2$, it follows the existence of a degree $2$ morphism $C\to\pp^1$,  contrary to the hypothesis that $C$ is nonhyperelliptic. Therefore $h^0(\ww(-P))>h^0(\ww(-2P))$. Thus there exists $y_1\in H^0(\ww(-P))\subset H^0(\ww)$ such that $v(y_1)=1$ and (ii) is proved.

Assume $C$ is hyperelliptic and $g\geq 2$, then there is $a\in k(C)$ such that $\text{div}_0(a)=P+Q$ and a differential $y$ such that $H^0(\ww)=\langle y,ay,\ldots,a^{g-1}y\rangle$. Again, since $\ww$ is generated by $H^0(\ww)$ and since $v_P(a)\geq 1$ one has $\wwp=\op y$ which forces $v_P(y)=0$. Setting $y_1:=by$ then $v_P(y_1)=1$ or $2$ depending on if $P$ differs from $Q$ or not. This proves (iii).
\end{proof}

\begin{lemma}
\label{lemchi}
If $g\geq 2$ then $h^1(\ww^n)=0$ for $n\geq 2$.
\end{lemma}

\begin{proof}
We first claim $\deg(\ww^n)>\deg(\ww)$ for every $n\geq 2$. In fact, set $\eta_P=\dim(\ww_P/\op x)$ and $\eta=\sum_{P\in C}\eta_P$. Then, focusing just on the contributions of Gorenstein points to $\deg(\ww)$ we get
\begin{align*}
\deg(\ww^n)-\deg(\ww)&\geq (n-1)(2g-2-\eta)\\
                     &=(n-1)((g-2)+(g-\eta)).
\end{align*}
But $g>\eta$ unless $g=0$ owing to the proof of \cite[Prp.\,5.2]{KM}. This proves the claim. Therefore $\chi(\ww^n)>\chi(\ww)=g-1$ and the result follows from \cite[Lem.\,5.1.(1)]{KM}.
\end{proof}
  
\begin{lemma}
\label{lemngp}
Let $P\in C$ be a unibranch non-Gorenstein point. Let $\ct$ be the curve obtained from $C$ resolving $P$ and $\gtl$ its genus. Let $\varphi :\ct\to C$ be the natural morphism and $\wwt:=\varphi_*(\ww_{\ct})$. If $\widetilde{g}\geq 2$, then
$$
\text{\emph{Sym}}^n H^0(\ww)\longrightarrow H^0(\ww^n)/H^0(\wwt^n)
$$
is surjective for every $n\geq 1$.
\end{lemma}

\begin{proof}
Form the long exact sequence
$$
0\to H^0(\wwt^n)\to H^0(\ww^n)\stackrel{u}\to H^0(\ww^n/\wwt^n)\to H^1(\wwt^n)\to H^1(\ww^n).
$$
Now $h^1(\ww)=h^1(\wwt)=1$ and $g>\gtl\geq 2$ what implies $h^1(\ww^n)=h^1(\wwt^n)=0$ for $n\geq 2$ due to Lemma \ref{lemchi}. Then $u$ is surjective and so $H^0(\ww^n)/H^0(\wwt^n)=H^0(\ww^n/\wwt^n)=\wwp^n/\wwtp^n$. So we need showing $\text{Sym}^nH^0(\ww)\to\wwp^n/\wwtp^n$ is surjective. 

By construction $\wwtp =\overline{\ww}_P$, and $\overline{\ww}_P=\mathcal{C}_Px$ owing to \cite[Lem.\,2.8]{KM}. This implies $\vtxp=\cp$. Now $\wwp^n/\wwtp^n=\vxp^n/\vtxp^n=\vxp^n/\cp^n$. Besides, $\text{Sym}\ H^0(\ww)\to W_x^n$ is surjective. So it suffices showing $W_x^n\to\vxp^n/\cp^n$ is surjective for $n\geq 1$.

Let $\pt$ be the point of $\ct$ which lies over $P$. Since $\gtl\geq 2$, Lemma \ref{lemvpw} implies the existence of $y_0,y_1\in H^0(\wwt)$ such that $v_{\pt}(y_0)=0$ and $v_{\pt}(y_1)=1$ or $2$. But $H^0(\wwt)\subset H^0(\ww)$ and $v_{\pt}=v_{P}$ as functions by the very definition. So there are $y_0,y_1\in H^0(\ww)$ such that $v_{P}(y_0)=0$ and $v_{P}(y_1)=1$ or $2$. Lemma \ref{lemtch}.(iii) implies $W_x^n\to\vxp^n/\cp^n$ is surjective for $n\geq 1$.
\end{proof}

\begin{lemma}
\label{lemhyp}
Let $P\in C$ be a unibranch point of multiplicity at least $3$. Let $\ct$ be the curve obtained from $C$ resolving $P$ and $\gtl$ its genus. Let $\varphi :\ct\to C$ be the natural morphism and $\wwt:=\varphi_*(\ww_{\ct})$. If $\ct$ is hyperelliptic and $\widetilde{g}\geq 2$, then 
$$
\text{\emph{Sym}}^n H^0(\ww)\rightarrow H^0(\wwt^n)
$$
is surjective for every $n\geq 1$.
\end{lemma}

\begin{proof}
Since $\ct$ is hyperelliptic and $\gtl\geq 2$, there is $a\in k(C)$ such that $\deg(\text{div}_0(a))=2$ and 
\begin{equation}
\label{equbas}
H^0(\wwt )=\langle y,ay ,\ldots ,a^{\gtl-1}y\rangle
\end{equation} 
for a certain $y\in H^0(\wwt )$. On the one hand, $\dim(H^0(\wwt)^n)= n(\gtl-1)+1$ due to (\ref{equbas}). On the other hand, if $n\geq 2$ then $\dim(H^0(\wwt^n))=h^0(\wwt ^n)=(2n-1)(\gtl-1)$ owing to Riemann-Roch and Lemma \ref{lemchi}. So, if $n\geq 2$, there are elements in $H^0(\wwt ^n)$ which are not in $H^0(\wwt )^n$. We will show these elements are in $H^0(\ww)^n$ by proving that $\widetilde{W}_{y,n}\subset W_y^n$. 

Let $\pt\in\ct$ be the point which lies over $P$ and let $\psi:=(1:a):\ct\to \pp^1$ be the degree $2$ morphism. We will consider two cases.  

\noindent\emph{Case 1.} $\pt$ is a ramification point of $\psi$.

Let $\pb\in\overline{C}$ be the point which lies over $\pt$ and $P$. In this case we may take 
\begin{equation}
\label{equpol}
\text{div}(a)=2\pb -\overline{Q}-\overline{R}
\end{equation}
where $\overline{Q}, \overline{R}$ lie over nonsingular points, say $Q,R$, of $C$. 

Since the multiplicity of $P$ is at least $3$, one might consider the partial normalization map $\phi:C^*\to C$ where $C^*$ is obtained from $C$ replacing $P$ by a point $\pbb$ of multiplicity $3$ for which the maximal ideal agrees with the conductor. In other words, if $\obr:=\phi_*(\oo_{C^*})$ then $\obr_P=t_P^3\obp$. Set $\wwb:=\phi_*(\ww_{C^*})$. If $h^1(\obr(Q+R))=h^1(\obr(Q))$ then $h^0(\obr(Q+R))>h^0(\obr(Q))\geq h^0(\obr)=1$ because $\deg(\obr(Q+R))=\deg(\obr(Q))+1$. Now if $Q^*$ and $R^*$ are the points of $C^*$ lying over $Q$ and $R$ then $\deg_{C^*}(\oo_{C^*}(Q^*+R^*))=2$ and $h^0(\oo_{C^*}(Q^*+R^*))=h^0(\obr(Q+R))\geq 2$. It follows that $C^*$ is hyperelliptic which cannot happen because it has a point $\pbb$ with multiplicity $3$ (actually $\pbb$ is even non-Gorenstein). 

Therefore $h^1(\obr(Q+R))<h^1(\obr(Q))$, that is, $H^0(\wwb(-Q-R))\subsetneqq H^0(\wwb(-Q))$. By construction $\ww^*_Q=\wwt_Q$ and $\ww^*_R=\wwt_R$. Besides, (\ref{equbas}) and (\ref{equpol}) implie $\wwt_Q=\oo_Q\, a^{\gtl-1}y$ and $\wwt_R=\oo_R\, a^{\gtl-1}y$. Then there exists $b\in k(C)$ such that $by \in H^0(\ww ^*)$ and with $v_{Q}(b)=-\gtl+2$ and $v_{R}(b)=-\gtl+1$. 

Now set
$$
B:=\{1,a,\ldots ,a^{n(\gtl-1)},ba^2,\ldots ,ba^{(n-1)(\gtl-1)}\}.
$$ 
We claim $\widetilde{W}_{y,n}=\langle B\rangle\subset W_{y}^n$. In fact, first, $B$ is a $k$-linearly independent set as one can easily see computing $v_Q$ and $v_R$ of its elements. Second $|B|=(2n-1)(\hat{g}-1)=h^0(\wwt^n)=\dim(\widetilde{W}_{y,n})$. Third, $B\subset\widetilde{W}_{y,n}$. Indeed, the $a^i$'s are clearly in $W_{y,n}$ and the $ba^i$'s are in $\widetilde{V}_{y,S}^n$ for $n\geq 2$ and $S\neq P$ because $\wwt$ and $\ww^*$ agree outside $P$, $by\in H^0(\ww^*)$ and so $b\in \widetilde{V}_{y,S}$ for $S\neq P$. Let us show the $ba^i$'s are in $\widetilde{V}_{y,P}$. From (\ref{equbas}) and (\ref{equpol}) we have $\wwtp=\obp y$, in particular $v_P(y)=0$ as seen in the proof of Lemma \ref{lemvpw}. Now $by\in H^0(\ww ^*)$ and so, according to (\ref{equvwp}), its pole at $P$ is at most $\alpha _{P^*}=3$. Hence $v_P(b)=v_P(by)\geq -3$. Therefore $v_P(ba^2)=v_P(a^2)+v_P(b)\geq 4-3>0$. So $ba^2\in\obp=\widetilde{V}_{y,P}$. This implies all the $ba^i$'s are in $\widetilde{V}_{y,P}$. So we have proved $\widetilde{W}_{y,n}=\langle B\rangle$. Since $H^0(\ww ^*)\subset H^0(\ww )$ it follows that $b\in W_{y}$ and since $H^0(\wwt)\subset H^0(\ww )$ it follows that the $a^i$'s are in $W_y$. So $B\subset W_{y}^n$ as desired and we are done with first case. 

\noindent\emph{Case 2.} $\pt$ is not a ramification point of $\psi$

Now we consider $C^*$ the curve obtained from $C$ replacing $P$ by a point $P^*$ of multiplicity $2$ for which the maximal ideal agrees with the conductor. We claim $C^*$ is nonhyperlliptic. Indeed, otherwise we would be able to take $a\in\oo _{C^*,P^*}$. But if so, subtracting $a$ by a suitable constant we may suppose it is not a unit in $\oo _{C^*,P^*}$ and so $v_{P}(a)=2$ but this contradicts the fact that $\pt$ is not a ramification point of $\psi$. Actually, it is easy to see that $|\oo _{C^*}\langle 1,a\rangle |$ is a $g_{3}^1$ on $C^*$ with a non-removable base point $P^*$ and so $C^*$ is trigonal. 

Then we proceed verbatim as in Case 1 up to the following change: we take $\text{div}(a)=\pb +\pb _1 -\overline{Q}-\overline{R}$ with  $\pb _1\neq\pb$ and so $v_P(a)=1$; we have $\alpha_{P^*}=2$ and hence $v_P(b)=v_P(by)\geq -2$; therefore $v_P(ba^2)=v_P(a^2)+v_P(b)\geq 2-2=0$. 
\end{proof} 

\begin{thm}
\label{thmmon}
Let $C$ be a nonhyperelliptic curve whose non-Gorenstein points are unibranch. Then the homomorphisms
$$
\text{\emph{Sym}}^{n}\,H^{0}(C,\ww )\longrightarrow H^{0}(C,\ww ^{n})
$$
are surjective.
\end{thm}

\begin{proof}
We proceed by induction on the number of non-Gorenstein points of $C$. Assume first there is only one non-Gorenstein point $P\in C$. Let $\ct$ be the curve obtained from $C$ resolving $P$ and $\gtl$ its genus. Let $\varphi :\ct\to C$ be the natural morphism and $\wwt:=\varphi_*(\ww_{\ct})$. 

If $\gtl=0$ then $h^0(\wwt^n)=0$ for $n\geq 1$ and $h^1(\ww^n)=0$ for $n\geq 2$ due to Lemma \ref{lemchi} because $g\geq 2$ since $C$ is non-Gorenstein. Besides, from Lemma \ref{lemtch}.(i) we have $\dim(W_x^n)\geq\dim(V_{x,P}^n/t_P^{-(2n-1)}\cp^n)$ and we have already seen $\widetilde{V}_{x,P}=\cp$. These statements lead us to
\begin{align*}
\dim(H^0(\ww^n))&=h^0(\ww ^n)=h^0(\ww ^n)-h^0(\wwt ^n) \\ 
                &=\deg(\ww ^n)-\deg(\wwt ^n)-h^1(\wwt ^n) \\
                &=\dim (\wwp ^n/\wwtp ^n)-(2n-1) \\ 
                &=\dim(V_{x,P}^n/\widetilde{V}_{x,P}^n)-(2n-1) \\
                &=\dim(V_{x,P}^n/\cp^n)-(2n-1) \\
                &=\dim(V_{x,P}^n/t_P^{-(2n-1)}\cp^n) \\
                &\leq\dim(W_x^n)=\dim(H^0(\ww)^n).
\end{align*}
Therefore the theorem holds if $C$ has just one non-Gorenstein point and $\gtl=0$.

If $\gtl=1$ then $\ww_{\ct}\cong\oo_{\ct}$. In particular, $H^0(\wwt)^n=H^0(\wwt^n)$. Besides, if $\gtl=1$ then there exits $y_0\in H^0(\wwt)\subset H^0(\ww)$ such that $v_P(y_0)=0$ owing to Lemma \ref{lemvpw}.(i). Hence Lemma \ref{lemtch}.(ii) implies $W_x^n\to V{x,P}^n/t^{-1}\cp^n$ is surjective. Since $\widetilde{V}_{x,P}=\cp$ and $\widetilde{W}_x=W_x\cap\widetilde{V}_{x,P}$ it follows that $W_x^n/\widetilde{W}_x^n\to V{x,P}^n/t_P^{-1}\cp^n$ remains surjective. These statements lead us to
\begin{align*}
\dim(H^0(\ww^n)/H^0(\wwt)^n)&=\dim(H^0(\ww^n)/H^0(\wwt^n))  \\
                            &=h^0(\ww^n)-h^0(\wwt^n) \\
                            &=\deg(\ww^n)-\deg(\wwt^n)-h^1(\wwt^n) \\
                            &=\dim(\wwp^n/\wwtp^n)-1 \\
                            &=\dim(V_{x,P}^n/\widetilde{V}_{x,P}^n)-1 \\
                            &=\dim(V_{x,P}^n/\cp^n)-1 \\
                            &=\dim(V_{x,P}^n/t_P^{-1}\cp^n) \\
                            &\leq \dim(W_x^n/\widetilde{W}_x^n) \\
                            &=\dim(H^0(\ww)^n/H^0(\wwt)^n). 
\end{align*}
Therefore the theorem holds if $C$ has just one non-Gorenstein point and $\gtl=1$.

Now $P$ has multiplicity at least $3$ for it is non-Gorenstein and $\ct$ is Gorenstein because $P$ is the only non-Gorenstein point of $C$. Thus if $\gtl\geq 2$ either $\ct$ is hyperelliptic and so $\text{Sym}^n\,H^{0}(\ww)\to H^{0}(\wwt^{n})$ is surjective owing to Lemma \ref{lemhyp} or else $\ct$ is nonhyperelliptic Gorenstein and so $\text{Sym}^n\,H^{0}(\wwt )\to H^{0}(\wwt^{n})$ is surjective owing to Theorem \ref{thmext}. Besides, $\text{Sym}^n\,H^{0}(\ww)\to H^0(\ww^n)/H^{0}(\wwt^{n})$ is surjective due to Lemma \ref{lemngp}. Therefore the theorem holds if $C$ has just one non-Gorenstein point.

If $C$ has many non-Gorenstein points, define $\ct$ resolving just one of them. Then $\text{Sym}^n\,H^{0}(\wwt )\twoheadrightarrow H^{0}(\wwt^{n})$ by induction and $\text{Sym}^n\,H^{0}(\ww)\twoheadrightarrow H^0(\ww^n)/H^{0}(\wwt^{n})$ owing to Lemma \ref{lemngp}. We are done.
\end{proof}

\begin{rem}
\emph{Though it was not our aim, one can combine Theorem \ref{thmext} with Theorem \ref{thmmon} to get the following result: if $C$ is a nonhyperelliptic curve with at most one multibranch non-Gorenstein point which is also almost Gorenstein, then $\text{Sym}^n\,H^{0}(\ww)\to H^{0}(\ww^{n})$ is surjective for $n\geq 1$. The proof is basically the same of Theorem \ref{thmmon} just using Theorem \ref{thmext} in its full generality.} 
\end{rem}

\end{document}